\documentclass[a4paper]{article}
\usepackage[utf8]{inputenc}  
\usepackage[T1]{fontenc}

\usepackage{lmodern} \normalfont 
\DeclareFontShape{T1}{lmr}{bx}{sc} { <-> ssub * cmr/bx/sc }{}

\usepackage{amsthm}
\usepackage{amssymb}
\usepackage{amsmath}
\usepackage{mathrsfs}
\usepackage{dsfont}
\usepackage{yhmath}
\usepackage{listings}

\usepackage{graphicx}

\usepackage{pgf,tikz}
\usetikzlibrary{arrows}
\usetikzlibrary{patterns}

\usepackage[all]{xy}

\usepackage{subcaption}
\usepackage{geometry}
\usepackage{ifthen}

\usepackage[colorlinks=true,linkcolor=blue,citecolor=green]{hyperref}

\newcommand{\N}{\mathbb{N}}
\newcommand{\Z}{\mathbb{Z}}
\newcommand{\Q}{\mathbb{Q}}
\newcommand{\C}{\mathbb{C}}

\renewcommand{\P}{\mathbb{P}}

\newcommand{\E}{{\cal E}}
\newcommand{\pic}{\mathrm{Pic}}

\def\lint{[\hskip -1.5pt[}
\def\rint{]\hskip -1.5pt]}

\def\bsys{\left\{\begin{array}}
\def\esys{\end{array}\right.}

\newcommand{\eps}{\varepsilon}

\renewcommand{\emptyset}{\varnothing}

\theoremstyle{plain}
\newtheorem{theorem}{Theorem}[section]
\newtheorem{lemma}[theorem]{Lemma}
\newtheorem{proposition}[theorem]{Proposition}
\newtheorem{corollary}[theorem]{Corollary}

\theoremstyle{definition}
\newtheorem{definition}[theorem]{Definition}
\newtheorem{notation}[theorem]{Notation}

\newtheorem{remark}[theorem]{Remark}
\newtheorem{example}[theorem]{Example}
\newtheorem{nexample}[theorem]{Non-exemple}

\newcommand{\lem}[2][None]
{\begin{lemma}\ifthenelse{\equal{#1}{None}}{}{\label{#1}}
#2
\end{lemma}}

\newcommand{\prop}[2][None]
{\begin{proposition}\ifthenelse{\equal{#1}{None}}{}{\label{#1}}
#2
\end{proposition}}

\newcommand{\cor}[2][None]
{\begin{corollary}\ifthenelse{\equal{#1}{None}}{}{\label{#1}}
#2
\end{corollary}}

\newcommand{\theo}[2][None]
{\begin{theorem}\ifthenelse{\equal{#1}{None}}{}{\label{#1}}
#2
\end{theorem}}

\newcommand{\defi}[2][None]
{\begin{definition}\ifthenelse{\equal{#1}{None}}{}{\label{#1}}
#2
\end{definition}}

\newcommand{\nota}[2][None]
{\begin{notation}\ifthenelse{\equal{#1}{None}}{}{\label{#1}}
#2
\end{notation}}

\newcommand{\rem}[2][None]
{\begin{remark}\ifthenelse{\equal{#1}{None}}{}{\label{#1}}
#2
\end{remark}}

\newcommand{\expl}[2][None]
{\begin{example}\ifthenelse{\equal{#1}{None}}{}{\label{#1}}
#2
\end{example}}

\newcommand{\nexpl}[2][None]
{\begin{nexample}\ifthenelse{\equal{#1}{None}}{}{\label{#1}}
#2
\end{nexample}}

\newcommand{\demo}[1]{\begin{proof} #1
\end{proof}}

\newcommand{\demode}[2]{\begin{proof}[Proof of #1] #2
\end{proof}}

\title{Positivity of the cotangent sheaf of singular Calabi-Yau varieties}

\author{Cécile Gachet\footnote{Université Côte d’Azur, CNRS, LJAD, France}}

\begin{document}

\maketitle

\begin{abstract}
We prove that the tangent and the reflexivized cotangent sheaves of any normal projective klt Calabi-Yau or irreducible holomorphic symplectic variety are not pseudoeffective, generalizing results of A. Höring and T. Peternell \cite{HoeringPeternell}. We provide examples of Calabi-Yau varieties of small dimension with singularities in codimension~2.
\end{abstract}


\section{Introduction}

Complex algebraic varieties with trivial canonical class are of great importance in birational geometry. Indeed, they appear naturally as possible minimal models in the Minimal Model Program (MMP), and come in quite diverse geometrical families. Since higher-dimensional MMP generally gives rise to singular minimal models, understanding singular projective varieties with trivial canonical class is particularly relevant. 
Recently, three papers \cite[Thm.1.5]{HoeringPeternell}, \cite{GGK}, \cite{Druel} achieved a singular decomposition result for these varieties :

\theo[BBDec]{Let $X$ be a normal projective variety with klt singularities, with $K_X$ numerically trivial. Then there exists a normal projective variety $\tilde{X}$ with at most canonical singularities, which comes with a quasiétale finite cover $f:\tilde{X}\to X$ and decomposes as a product:
$$\tilde{X}\cong A\times\prod_{i} Y_i\times\prod_{j} Z_j,$$
where $A$ is a smooth abelian variety, the $Y_i$ are singular Calabi-Yau varieties and the $Z_j$ are singular irreducible holomorphic symplectic (IHS) varieties, as defined in Section \ref{sec:thmnotpseff}.}

May it seem an expected generalization of the smooth Beauville-\linebreak Bogomolov decomposition result \cite{Beauv81}, \cite{Beauv}, this theorem however relies on serious results from each paper: \cite{GGK} introduces algebraic holonomy and studies infinitesimal decompositions of the tangent sheaf $\mathcal{T}_X$; \cite{Druel} deals with the abelian part in the infinitesimal decomposition through a positive characteristics argument, and proves an integrability criterion for the remaining subsheaves of $\mathcal{T}_X$; \cite{HoeringPeternell} establishes positivity results which add up to Druel's criterion to finish the proof. This proof was notably simplified by \cite{Campana}, shortcutting the positive characterictics argument. Furthermore, the two recent papers \cite{ClaudonGrafGuenancia}, \cite{BGL} extend this decomposition result to the singular Kähler case by subtle algebraic approximation considerations.

\medskip

Interestingly enough, the singular decomposition for a klt variety $X$ may not be the same as the singular decomposition of its terminalisation. The typical example is a singular Kummer surface, which resolves by 16 blow-ups into a smooth K3 surface, but has the Beauville-Bogomolov type of an abelian surface. Other such intruiguing examples are given in \cite[Sect.14]{GGK}. Compatibility of the singular Beauville-Bogomolov decomposition with terminalisation nevertheless holds for {\it some} klt varieties with trivial canonical class \cite[Lem.4.6]{Druel}. This license to terminalise is essential in the current proof of \cite[Thm.1.5]{HoeringPeternell}, as it involves positivity results \cite[Thm.1.1]{HoeringPeternell} for klt varieties which are smooth in codimension 2: any klt variety is not, but its terminalisation surely is.

Since these positivity results have a wider scope than the mere proof of the singular decomposition theorem, it is worth extending them to normal projective klt varieties. Our main theorem is:

\theo[CYNotPseff]{Let $X$ be a normal projective variety with klt singularities and numerically trivial $K_X$. 
If its tangent or reflexivized cotangent sheaf is pseudoeffective, then there is a quasiétale finite cover $\tilde{X}\to X$ such that $q(\tilde{X})\ne 0$. Equivalently, the singular Beauville-Bogomolov decomposition of $X$ has an abelian factor of positive dimension.

In particular, if $X$ is a singular Calabi-Yau or IHS variety in the sense of Def.\ref{CY}, then neither $\mathcal{T}_X$ nor its dual $\Omega^{[1]}_X$ is pseudoeffective.}

Importantly enough, this theorem does not boil down to \cite[Thm.1.6]{HoeringPeternell} on a terminalisation of $X$; we inevitably have to deal with codimension 2 quotient singularities on $X$. In this perspective, we resort to the theory of orbifold Chern classes. It has been developped in the late eighties in connection to the abundance problem for threefolds \cite{FlipsNA}, and we will extensively use some of its most recent developments, {\it inter alia} \cite{LuTaji}, \cite{GKPT}, \cite{GKPPreprint}.

\medskip

Let us present a brief outline of the proof, say for a variety $X$ with pseudoeffective tangent sheaf.

The fact that $\mathcal{T}_X$ is pseudoeffective pullbacks and restricts to one factor in the Beauville-Bogomolov decomposition of $X$. Supposing by contradiction that $X$ has no abelian part, we can reduce to a Calabi-Yau or IHS factor $Z$ such that $\mathcal{T}_Z$ is pseudoeffective. The work of \cite{GGK} also establishes that $\mathcal{T}_Z$ and all its symmetric powers are stable of slope zero with respect to any polarisation $H$. Finally, since $Z$ is not abelian, its orbifold second Chern class satisfies $\hat{c}_2(\mathcal{T}_Z)\cdot H^{\mathrm{dim}\, X - 2}\ne 0$. This contradicts the following generalization of \cite[Thm.1.1]{HoeringPeternell}:

\theo[main]{Let $X$ be a normal projective variety with klt singularities of dimension $n$, $H$ a $\Q$-Cartier ample divisor on $X$. Consider $\E$ a reflexive sheaf on $X$ such that:
\begin{itemize}
\item $\hat{c_1}(\E)\cdot H^{n-1} = 0$;
\item the sheaves $\E$ and $S^{[l]}\E$, for some $l\ge 6$, are $H$-stable;
\item $\E$ is pseudoeffective.
\end{itemize}

\noindent Then $\hat{c_1}(\E)^2\cdot H^{n-2} = \hat{c_2}(\E)\cdot H^{n-2} = 0.$

Moreover, there is a finite Galois covering $\nu:\tilde{X}\to X$, étale in codimension 1, such that $\nu^{[*]}\E$ is locally-free, has a numerically trivial determinant, and is $\mathrm{Gal}(\tilde{X}/X)$-equivariantly flat on $\tilde{X}$, {\it ie} comes from a $\mathrm{Gal}(\tilde{X}/X)$-equivariant representation of $\pi_1(\tilde{X})$. In particular, $\nu^{[*]}\E$ is numerically flat, and, as symmetric multilinear forms on $NS(X)$:
$$c_1(\nu^{[*]}\E)\equiv 0,\quad c_2(\nu^{[*]}\E)\equiv 0.$$}

The hard part here is the first assertion on the vanishing of orbifold Chern classes, the rest follows from \cite{LuTaji}. 

\medskip

In Section \ref{sec:foundations}, we recall and prove basics to reduce the proof of Thm.\ref{main} to working on a normal projective klt surface $S$. A crucial ingredient is that orbifold Chern classes behave well under certain restrictions \cite[Prop.3.11]{GKPT}.
In Section \ref{sec:surface}, we introduce an {\it unfolding} $p:\hat{S}\to S$, obtained by gluing together local finite Galois quasiétale resolutions of the singularities of $S$. The surface $\hat{S}$ may be as singular as $S$; importantly enough though, any reflexive sheaf $\E$ on $S$ reflexively pulls back to a locally-free sheaf $\hat{\E}$ on $\hat{S}$. We investigate the relationship of $\E$ and $\hat{\E}$. The key of the proof of Thm.\ref{main} is then to establish the nefness of $\hat{\E}$, which yields the Chern classes vanishing for $\hat{\E}$, hence for $\E$. Note that $\E$ may very well not be nef itself: see Remark \ref{rem:interesting}.

As a conclusive remark, note that investigating pseudoeffectivity of the tangent and reflexivized cotangent sheaves of a variety with trivial canonical class requires knowledge of its singular Beauville-Bogomolov decomposition. To that extent, Thm.\ref{CYNotPseff} cannot be used on an explicit variety before knowing a bare minimum about its geometry. In Section \ref{sec:threeexpl}, we exhibit 2409 Calabi-Yau threefolds with singularities in codimension 2, among the 7555 wellformed quasismooth hypersurfaces of trivial canonical sheaf in weighted projective 4-dimensional spaces classified by \cite{KreuzerSkarke}. These examples stay out of the range of the earlier pseudoeffectivity result of \cite[Thm.1.6]{HoeringPeternell}, but are covered by our Thm.\ref{CYNotPseff}.

\subsection*{Acknowledgments} 
I heartily thank my advisor A. Höring for suggesting me to study this subject, for fruitful discussions and for his careful reading of several versions of this paper. My gratitude also goes to the anonymous referee, for their noteworthy remarks on how to make this paper more concise and readable. 

\section{Notation and basic facts}\label{sec:foundations}

\subsubsection*{Finite morphisms} We will deal with various types of finite maps.

\defi{Unless otherwise stated, all finite morphisms we speak about are surjective; we may well refer to them as finite coverings, without saying anything about how étale they are.
\noindent We refer say that a finite morphism is {\it quasiétale} if it is étale in codimension 1.
\noindent Following \cite{GKPFund}, we call a finite morphism of normal varieties $Y\to X$ {\it Galois} if it is the quotient map of $Y$ by a finite group action.}

\subsubsection*{Reflexive sheaves} Let $\E$ be a reflexive sheaf on a variety $X$. Recall the reflexivization functor $\mathcal{F}\mapsto \mathcal{F}^{**}$ enables to perform general algebraic operations in the category of reflexive sheaves. 

Notably, we will denote by:
\begin{itemize}
\item $S^{[l]}\E$ the reflexivization of the $l$-th symmetric power of $\E$ (for $l\in\N$), 
\item $\nu^{[*]}\E$ the reflexivization of the pullback of $\E$ (for $\nu: Y\to X$ a morphism).
\end{itemize}

Note that, by \cite[Prop.1.6]{HartRefl}, reflexive sheaves are normal: if $\E$ is a reflexive sheaf on $X$, then for all open sets $V\subset U\subset X$ such that $\mathrm{codim}_U U\setminus V\ge 2$, the restriction map $\E(U)\to\E(V)$ is an isomorphism.
In particular, a morphism between two reflexive sheaves which is an isomorphism when restricted to a {\it big open set} (ie an open set whose complementary has codimension at least 2) is a global isomorphism.

\lem[ReflPB]{Let $p:X\to Y$ be a finite morphism between normal projective varieties. The functor $p^{[*]}$ from the category of reflexive sheaves on $Y$ to that of reflexive sheaves on $X$ is left-exact.}

\demo{Let $0\to \E\to \mathcal{F} \to \mathcal{G}$ be an exact sequence of reflexive sheaves on $Y$, and denote by $Z\subset Y$ a closed subscheme of codimension at least 2 such that our reflexive sheaves are locally-free on $Y\setminus Z\subset Y_{\mathrm{reg}}$. Reflexive pullback a priori only gives morphisms
$$p^{[*]}\E\to p^{[*]}\mathcal{F} \to p^{[*]}\mathcal{G},$$ whose composition is zero.
By \cite[Lem.31.12.7]{StackEx}, the kernel $K$ of the morphism $p^{[*]}\mathcal{F}\to p^{[*]}\mathcal{G}$ is reflexive. There is a natural morphism from $p^{[*]}\E$ to the kernel $K$, which restricts to an isomorphism over $X\setminus p^{-1}(Z)$. As both sheaves are reflexive and $p^{-1}(Z)$ has codimension 2, they are isomorphic over all $X$.}

\subsection{Nefness and pseudoeffectivity}


Let us recall that a coherent sheaf $\E$ on a normal variety $X$ has a projectivization $\P(\E)$ with a canonical, so-called tautological, line bundle $\zeta$ on it and a natural morphism $p:\P(\E)\to X$ with a natural sheaf quotient map: $p^*\E\twoheadrightarrow\zeta$. An account on this set-up is given in \cite[Chapt.4]{EGA2}.
We simply recall the universal property of this construction: for any scheme $q:C\to X$, to give an $X$-morphism $\nu: C\to \P(\E)$ is equivalent to giving a line bundle $L$ over $C$ together with a sheaf surjection $q^*\E\twoheadrightarrow L$.

Projectivizations are standardly used for generalizing positivity notions of line bundles to coherent sheaves, as follows.

\defi{Let $\E$ be a coherent sheaf on a normal variety $X$. It is called {\it nef} if the tautological line bundle $\zeta$ on $\P(\E)$ is nef.}

\rem{This coincides with \cite[Def.6.1.1]{Laz2} when the sheaf $\E$ is locally-free.
Note that for a torsion-free coherent sheaf $\E$, the scheme $\P(\E)$ may well have several irreducible components. Somehow, several of these components may be relevant for studying the positivity of $\E$: not only the mere one which is dominant onto $X$, but also components which may be contracted to a non-zero dimensional locus of $X$. Such components don't exist for a reflexive sheaf on a normal projective surface: so in this case, nefness is easier to study.}

\prop[moreonnef]{We have the following properties:
\begin{itemize}
\item if $Y\subset X$ is a normal subvariety, and $\E$ is a nef coherent sheaf on $X$, then $\E|_Y$ is nef;
\item conversely, nefness of a coherent sheaf $\E$ is enough to be checked on all curves of $\P(\E)$;
\item if $f:Y\to X$ is a finite dominant morphism of normal varieties and $\E$ is a coherent sheaf on $X$, $\E$ is nef if and only if $f^*\E$ is;
\item if $f:Y\to X$ is a proper birational morphism resolving the singularities of a normal variety $X$ and $\E$ is a coherent sheaf on $X$ such that $f^*\E$ is nef, then $\E$ is nef.
\end{itemize}

\noindent These are simple consequences of the universal property of $\P(\E)$ and of the fact \cite[4.1.3.1]{EGA2} that for a dominant morphism $f:Y\to X$ and a coherent sheaf $\E$ in $X$, 
$$\P(f^*\E)=\P(\E)\times_X Y.$$}

\rem[rem:interesting]{Interestingly enough, the reflexive pullback of a non-nef reflexive sheaf may be nef, as shows the following example. Let $X$ be a {\it singular Kummer surface}, ie the finite quotient of an abelian surface $A$ by the involution $i:a\mapsto -a$.
Since $p:A\to X$ is a finite quasiétale cover and $\mathcal{T}_X$ is locally-free on a big open set, the reflexive sheaves $p^{[*]}\mathcal{T}_X$ and $\mathcal{T}_A$ are the same. In particular, 
$$p^{[*]}\mathcal{T}_X=\mathcal{O}_A\oplus\mathcal{O}_A\mbox{ is nef.}$$

We are going to prove that $\mathcal{T}_X$ itself is not nef. We first compute it.

Recalling \cite[App.B]{GrebKebekusKovacsPeternell}, we consider the functor taking the invariant direct image of a $\Z_2$-coherent sheaf on $A$: $\E\mapsto p_*\E^{\Z_2}$. It sends reflexive sheaves to reflexive sheaves, so that the following equality, which is clear on the big open étale locus of $p$, extends to a global sheaf isomorphism:
\begin{align}\label{pbpf}
&\text{for any $\E$ reflexive sheaf on $X$, $p^{[*]}\E$}\\
\notag &\text{is naturally $\Z_2$-equivariant and }(p_*p^{[*]}\E)^{\Z_2}\cong \E.
\end{align}

Moreover, it is an exact functor.
Still from \cite[App.B]{GrebKebekusKovacsPeternell}, if $\E$ is a ${\Z_2}$-equivariant coherent sheaf on $A$, then the sheaf $p_*\E^{\Z_2}$ is a direct summand of $p_*\E$. Note that a given coherent sheaf $\E$ on $A$ may have several structures of $\Z_2$-equivariant object. 
For example, $\mathcal{O}_A$ comes with a natural and a reversed action, defined on an affine open set $U$ by:
$$f\in\mathcal{O}_A(U)\mapsto f\circ i\in\mathcal{O}_A(i(U)),$$
$$f\in\mathcal{O}_A(U)\mapsto -f\circ i\in\mathcal{O}_A(i(U)).$$

\noindent So $p_*\mathcal{O}_A$ is the direct sum of two reflexive sheaves of rank 1, $\mathcal{O}_X$ (the invariants by the natural action) and $\mathcal{F}$ (the invariants by the reversed action).

We know that $\mathcal{T}_X \cong p_*(\mathcal{O}_A\oplus \mathcal{O}_A)^{\Z_2}$ with Eq.\ref{pbpf}. We note that the natural $\Z_2$-equivariant structure on $\mathcal{T}_A$ acts diagonally, and reversely on each trivial summand. So $\mathcal{T}_X \cong\mathcal{F}\oplus\mathcal{F}$.

Let us finally check that $\mathcal{F}$ is not nef. We compute locally: let $V\subset X,U=p^{-1}(V)\subset A$ be affine open sets with local coordinates $(x,y)\in\C^2\simeq U$ so that $p|_U$ ramifies only at $(0,0)$. The quotient map $p:U\to V$ rewrites:
\begin{align*}
\C[u,v,w]/(uv-w^2)\cong\mathcal{O}_X(V) &\to \C[x,y]\cong\mathcal{O}_A(U) \\
u,v,w &\mapsto x^2,y^2,xy,
\end{align*}

\noindent so its image $\C[ x^2,y^2,xy]$ identifies with the local ring $\mathcal{O}_X(V)$.
Hence,
\begin{align*}\mathcal{F}(V) 
&\simeq \{f\in\C[ x,y]\mid\forall\, x,y,\, f(x,y)= -f(-x,-y)\}\\
&= x\,\C[ x^2,y^2,xy]\oplus y\,\C[ x^2,y^2,xy],\end{align*}

\noindent so that $\mathcal{F}^{\otimes 2}(V) \simeq u\,\mathcal{O}_X(V)\oplus v\,\mathcal{O}_X(V)\oplus w\,\mathcal{O}_X(V) = \mathcal{I}_{\mathrm{Sing}(X)}(V)$. This isomorphism is actually global: $$\mathcal{F}^{\otimes 2}\cong \mathcal{I}_{\mathrm{Sing}(X)}.$$ 
Ideal sheaves are not nef, so $\mathcal{F}^{\otimes 2}$ is not nef, so by \cite[Prop.2]{Kubota}, $\mathcal{F}$ is not nef. \qed}

\medskip


Pseudoeffectivity is standardly defined for locally-free sheaves through projectivisation too:

\defi[PsLocFr]{Le $\E$ be a locally-free sheaf on a normal projective variety $X$. It is considered {\it pseudoeffective} if it satisfies one of the following equivalent \cite[Lem.2.7]{Druel} conditions:
\begin{itemize}
\item the tautological line bundle on $\P(\E)$ is pseudoeffective;
\item there is an ample Cartier divisor $H$ on $X$ such that for all $c>0$, there are integers $i,j$ such that $i>cj>0$ and
$$h^0(X,\mathrm{Sym}^{i}\E\otimes \mathcal{O}_X(jH))\ne0.$$
\end{itemize}}

Generalizing this definition to any coherent sheaf is not obvious \cite{HPnew}. For reflexive sheaves, we use \cite[Def.2.1]{HoeringPeternell}: let $X$ be a normal projective variety and $H$ an ample Cartier divisor on $X$. A reflexive sheaf $\E$ on $X$ is said {\it pseudoeffective} if, for all $c>0$, there are numbers $i,j\in\N$ with $i>cj$ such that:
$$h^0(X, S^{[i]}(\E)\otimes\mathcal{O}_X(jH))\ne 0.$$

\expl{The sheaf $\mathcal{T}_X$ in Remark \ref{rem:interesting} is pseudoeffective, as $\mathcal{T}_X=\mathcal{F}\oplus \mathcal{F}$ with $S^{[2]}\mathcal{F}\cong\mathcal{O}_X$.}

\defi[PseffNaka]{Let $\E$ be a reflexive sheaf on a normal projective variety $X$. Denote by $\P'(\E)$ the normalization of the unique dominant component of $\P(\E)$ onto $X$. Let $P$ be a resolution of $\P'(\E)$, such that the birational morphism $r:P\to\P'(\E)$ over $X$ is an isomorphism precisely over the open locus $X_0\subset X_{\mathrm{reg}}$ where $\E$ is locally-free.

Denoting by $\pi$ the morphism $P\to\P(\E)$ and by $\mathcal{O}_P(1)$ the pullback of the tautological bundle of $\P(\E)$ by $\pi$, \cite[V.3.23]{Nakayama} asserts that one can choose (often not uniquely) an effective divisor $\Lambda$ supported in the exceptional locus of $r$ such that 
$$\zeta:=\mathcal{O}_P(1)\otimes\mathcal{O}_P(\Lambda)$$
satisfies $\pi_*\zeta^{\otimes m}\simeq S^{[m]}\E$ for all $m\in\N$. Such $\zeta$ is called a {\it tautological class} of $\E$.

As said in \cite[Lem.2.3]{HoeringPeternell}, with the same notations as previously, $\zeta$ is pseudoeffective on $P$ if and only if $\E$ is pseudoeffective as a reflexive sheaf.}

\prop[PseffRestr]{Let $X$ be a normal projective variety, $H$ an ample $\Q$-Cartier divisor, $\E$ a pseudoeffective reflexive sheaf on $X$. Then for $m$ big and divisible enough, for a general element $D\in |mH|$, the sheaf $\E|_D$ is reflexive and pseudoeffective.}

\demo{Let $U\subset X_{\mathrm{reg}}$ be a big open set on which $\E$ is locally-free. For $m$ big and divisible enough and for a general element $D$ in $|mH|$, $U\cap D$ is a big open set of $D$. By \cite[Prop.5.1.2]{GKPFund}, we can assume $D$ is a normal subvariety of $X$ and $\E|_D$ is reflexive.

Let us fix a $c>0$ and take $i,j$ integers such that $i> cj> 0$ and\linebreak $h^0(X,S^{[i]}(\E)\otimes\mathcal{O}_X(jH))> 0$. Up to taking a smaller $j$ (which may be negative if needed), we can assume that $$h^0(X,S^{[i]}(\E)\otimes\mathcal{O}_X((j-m)H))=0.$$
By normality of reflexive sheaves, 
\begin{align*}
h^0(D, S^{[i]}(\E|_D)\otimes\mathcal{O}_D(jH))
&= h^0(U\cap D, S^{i}(\E|_{U\cap D})\otimes\mathcal{O}_{U\cap D}(jH)) \\
&\ge h^0(U,S^{i}(\E|_{U})\otimes\mathcal{O}_U(jH)) \\
&\qquad - h^0(U,S^{i}(\E|_{U})\otimes\mathcal{O}_U((j-m)H))\\
&= h^0(X,S^{[i]}(\E)\otimes\mathcal{O}_X(jH)) \\
&\qquad - h^0(X,S^{[i]}(\E)\otimes\mathcal{O}_X((j-m)H))\\
& > 0,
\end{align*}
where the second equality comes from tensoring by $S^{i}(\E|_{U})\otimes\mathcal{O}_U(jH)$ and going to cohomology in the following exact sequence:
$$0\to \mathcal{O}_U(-mH)\to \mathcal{O}_U\to \mathcal{O}_{U\cap D}\to 0.\vspace{-2em}$$
}

We will use several times the following result \cite[Lem.3.15]{HPnew}:

\prop[PseffPullb]{Let $\E$ be a reflexive sheaf on a normal projective variety $X$, and $f:Y\to X$ be a finite dominant morphism of normal projective varieties. If $\E$ is pseudoeffective, then $f^{[*]}\E$ is.}

\defi{Let $D$ be a $\Q$-Cartier divisor on a normal projective variety $X$. We define its {\it stable base locus}:
$$\mathbb{B}(D):=\bigcap_{m\in M}\mathrm{Bs}(mD),$$
where $M\subset\N$ is the set of all $m$ such that $mD$ is Cartier and $\mathrm{Bs}(mD)$ is the base locus of the linear system $|mD|$.

\noindent We then define its {\it restricted base locus}:
$$B_-(D):=\bigcup_{n\in \N^*}\mathbb{B}\left(D+\frac{1}{n} A\right),$$
where $A$ is an arbitrary very ample divisor. Note that the union is strictly decreasing.}

Of course, a nef $\Q$-divisor has an empty restricted base locus. To that extent, the restricted base locus measures the non-nefness of a pseudoeffective line bundle. However, not all curves of a restricted base locus $B_-(D)$ must be $D$-non-positive, even in the simpler case where $D$ is a line bundle on a smooth surface and $B_-(D)$ is the negative part of its Zariski decomposition.

\subsection{Stability}

Let $\E$ be a torsion-free coherent sheaf on a normal projective variety $X$. For any ample $\Q$-Cartier divisor $H$ on $X$, and for some $m$ big and divisible enough, $n-1$ general members of $|mH|$ cut out a smooth curve $C$ on which $\E|_C$ is locally-free. So the usual notions of {\it slope stability} and {\it slope semistability} for $\E$ with respect to $H$ make good sense.

A generalization of a well-known Mehta-Ramanathan result says that stability behaves well under some well-chosen restrictions; we recall it as it is stated in \cite[Lem.2.11]{HoeringPeternell}:

\lem[Mehta]{Let $X$ be a normal projective variety of dimension $n$, and $H$ an ample Cartier divisor on $X$. Let $\E$ be a torsion-free coherent sheaf on $X$, that is stable with respect to $H$. Then there is $m_0$ such that, for all $m\ge m_0$, and for $D_1,\ldots,D_k$ general elements of $|mH|$ with $k\in\lint 1, n-1\rint$, if we denote by $Y$ the complete intersection $D_1\cap\ldots\cap D_k$, $\E|_Y$ is stable with respect to $H|_Y$.}

\rem{Note that the converse is clearly true.}

Stability a priori weakens through finite Galois reflexive pullbacks:

\lem[SemisPullb]{Let $p:Y\to X$ be a finite Galois cover of normal projective varieties of dimension $n$, $G$ its Galois group, $H$ an ample $\Q$-Cartier on $X$, $\E$ be a reflexive sheaf on $X$. Let $\mathcal{F}:=p^{[*]}\E$.
Then, if $\E$ is $H$-stable, $\mathcal{F}$ is $p^*H$-semistable.}

\demo{Suppose that $\E$ is $H$-stable. By Lemma \ref{Mehta}, on a smooth curve $C$ cut out by $n-1$ very general elements of the linear system defined by a suitable multiple of $H$, the now locally-free sheaf $\E|_C$ is still $H|_C$-stable.
In particular, \cite[Lem.6.4.12]{Laz2} applies; so the pullback sheaf $\mathcal{F}|_{p^{-1}(C)}$ is $p^*H|_C$-semistable. Hence, $\mathcal{F}$ is $H$-semistable.
}

Note that positivity and stability of a zero-slope locally-free sheaf are related by Miyaoka's result \cite{Miyaoka}, \cite[Prop.6.4.11]{Laz2}:

\prop[Miya]{Let $\E$ be a vector bundle on a smooth curve $C$. If $\E$ is semistable and $c_1(\E)=0$, then $\E$ is nef.}

More subtle than the mere stability of $\E$ is the stability of $\E$ and some of its symmetric powers:

\rem[BK]{We recall an interesting fact stated in \cite[Cor.6, following Rmk.]{BalajiKollar}. If $\E$ is a locally-free stable sheaf on a smooth projective variety $X$, then the following are equivalent:
\begin{itemize}
\item $S^r\E$ is stable for some $r\ge 6$ ;
\item $S^r\E$ is stable for any $r\ge 6$.
\end{itemize}

\noindent Whether or not the stability of all $S^{[l]}\E$ for $l\in\N$ could boil down to the stability of some $S^{[l]}\E$ for a finite amount of $l$'s remains an open question, when asked about a reflexive sheaf $\E$ on a smooth projective variety $X$ or about a locally-free sheaf $\E$ on a normal projective variety $X$.}

Nevertheless, this remark allows us to rewrite the key result \cite[Prop.1.3]{HoeringPeternell} in the following way:

\lem[Key]{Let $\E$ be a locally-free sheaf on a smooth curve $C$. Assume that $\E$ and $S^l\E$ for some $l\ge 6$ are stable and that $c_1(\E)=0$. Denoting by $\zeta$ the tautological bundle on $\P(\E)$, $\zeta$ is nef and satisfies:
$$\zeta^{\,\mathrm{dim}\, Z}\cdot Z > 0,$$
for any closed proper subvariety $Z\subset\P(\E)$.}

Despite that the reflexive pullback $p^{[*]}\E$ of a $H$-stable reflexive sheaf $\E$ by a finite dominant morphism $p$ is merely $p^*H$-semistable and a priori not stable (let alone his reflexive symmetric powers), the conclusive property of Lemma \ref{Key} is preserved by $p^{[*]}$:

\rem[QuasiBig]{Let $\E$ be a reflexive sheaf on a normal projective variety $X$, and $C\subset X$ a smooth curve such that $\E$ is locally-free in an analytical neighborhood of $C$, such that the tautological bundle $\zeta$ on $\P(\E|_C)$ is nef and such that it holds:
$$\zeta^{\mathrm{dim}\, Z}\cdot Z > 0,$$
for any closed proper subvariety $Z\subset\P(\E|_C)$.

Let $p:\hat{X}\to X$ be a finite dominant morphism, where $\hat{X}$ is a normal projective variety. Denote $\hat{C}:=p^{-1}(C)$, $\hat{\E}:=p^{[*]}\E$ and $\hat{\zeta}$ the tautological bundle of $\P(\hat{\E}|_{\hat{C}})$. If we have that $p^*(\E|_C)=\hat{\E}|_{\hat{C}}$, then the following diagram is Cartesian with tautological compatibility $\hat{\zeta}=q^*\zeta$:

\[
\xymatrix{
\P(\hat{\E}|_{\hat{C}}) \ar[d]^*{\hat{\pi}} \ar[r]^*{q} & \P(\E|_C) \ar[d]^*{\pi}\\
\hat{C} \ar[r]^*{p} & C
} 
\]

\noindent Hence, $\hat{\zeta}$ is nef and satisfies, for any closed proper subvariety $Z\subset\P(\hat{\E}|_{\hat{C}})$:
$$\hat{\zeta}^{\,\mathrm{dim}\, Z}\cdot  Z > 0.$$}

\rem{The case in which this remark will be relevant for us is when $X$ is a normal projective surface with an ample $\Q$-Cartier divisor $H$, $C$ is a smooth curve arising as a very general element of $|mH|$, for $m$ big and divisible enough, and $p:\hat{X}\to X$ is the morphism constructed in Section \ref{subsec:orbi}, so that $\hat{\E}$ is locally-free.
In this set-up, \cite[Prop.3.11.1]{GKPT} grants the additional assumption $p^*(\E|_C)=\hat{\E}|_{\hat{C}}$.}

\subsection{Orbifold Chern classes}\label{subsec:orbi}

Here we recall a standard construction for orbifold first and second Chern classes of a reflexive sheaf $\E$ on a normal projective variety $X$, whose singularities in codimension 2 are all quotient singularities. References for this matter include \cite{FlipsNA}, \cite{LuTaji}, \cite{GKPT}, \cite{GKPPreprint}. Note that normal projective klt varieties fall into this framework by the classical result \cite[Cor.1.14]{ReidCan}, \cite[Prop.9.4]{GrebKebekusKovacsPeternell}.

\medskip

Let $X$ be a normal projective variety, whose singularities in codimension 2 are all quotient singularities. Let $H$ be an ample $\Q$-Cartier divisor and $\E$ a reflexive sheaf on $X$. There are a normal quasiprojective subvariety $Y$ of $X$ with $\mathrm{codim}_X(X\setminus Y)\ge 3$ admitting an orbifold structure, a normal quasiprojective variety $\hat{Y}$, and a finite Galois morphism $p:\hat{Y}\to Y$ with Galois group $G$, such that $\hat{\E}:=p^{[*]}\E$ is a locally-free $G$-equivariant sheaf on $\hat{Y}$. 
Let us call the whole data $(X,Y,\hat{Y},p)$ an {\it unfolding} of $X$.

We can now define a first, respectively a squared first and a second orbifold Chern class of $\E$ as multilinear forms on $NS(X)^{n-1}$, respectively $NS(X)^{n-2}$ by:

\begin{align*}
\hat{c_1}(\E)\cdot H_1\cdots H_{n-1} 
&= \frac{1}{m^{n-1}\cdot |G|}c_1(\hat{\E})\cdot p^*(mH_1)\cdots p^*(mH_{n-1}),\\
\hat{c_1}^2(\E)\cdot H_1\cdots H_{n-2} 
&= \frac{1}{m^{n-2}\cdot |G|}c_1(\hat{\E})^2\cdot p^*(mH_1)\cdots p^*(mH_{n-2}),\\
\hat{c_2}(\E)\cdot H_1\cdots H_{n-2}
&= \frac{1}{m^{n-2}\cdot |G|}c_2(\hat{\E})\cdot p^*(mH_1)\cdots p^*(mH_{n-2}),
\end{align*}
where $H_1,\ldots,H_{n-1}$ are ample $\Q$-classes, and $m$ is big and divisible enough that general elements of $p^*(mH_1),\ldots,p^*(mH_{n-1})$ cut out a complete intersection smooth curve in $\hat{Y}$ and general elements of $p^*(mH_1),\ldots,p^*(mH_{n-2})$ a complete intersection normal surface in $\hat{Y}$.

As stated in \cite[Thm.3.13.2]{GKPT}, these orbifold Chern classes are compatible with general restrictions, and so is the unfolding construction \cite[Prop.3.11]{GKPT}.

\section{Restricting to a general surface}\label{sec:surface}

We prove the following proposition in Section \ref{subsec:techn}: 

\prop[OnSurface]{Let $S$ be a normal projective klt surface, and $H$ an ample $\Q$-Cartier divisor on $S$. Let $\E$ be a reflexive sheaf on $S$ such that:
\begin{itemize}
\item $\hat{c_1}(\E)\cdot H = 0$;
\item $\E$ and $S^{[l]}\E$, for some $l\ge 6$, are stable with respect to $H$;
\item $\E$ is pseudoeffective.
\end{itemize}

\noindent Then there is an unfolding $p:\hat{S}\to S$ as in Section \ref{subsec:orbi} on which the locally-free sheaf $\hat{\E}=p^{[*]}\E$ is nef.}

In Section \ref{subsec:surf}, we explain how this result implies the first part of Theorem \ref{main}, namely the vanishing of the squared first and second orbifold Chern classes.

\subsection{Consequences of Proposition \ref{OnSurface}}\label{subsec:surf}

We are going to combine Proposition \ref{OnSurface} with this lemma:
\lem[NefChern]{Let $S$ be a normal projective surface, $H$ an ample $\Q$-Cartier divisor on $S$ and $\E$ a locally-free sheaf on $S$. Assume that $\E$ is nef and $c_1(\E)\cdot H=0$. Then:
$$c_1(\E)^2=c_2(\E)=0.$$}

\demo{Let $\tilde{S}\overset{\varepsilon}{\to} S$ be the minimal resolution of $S$, $\tilde{H}=\varepsilon^*H$.
Writing $\tilde{\E}:=\varepsilon^*\E$, we get a nef locally-free sheaf on a smooth surface. The functoriality of Chern classes of locally-free sheaves by continuous pullbacks \cite[XI-Lem.1]{MilnorStashev} guarantees
$c_i(\tilde{\E})=\varepsilon^* c_i(\E)$ for $i=1,2$. In particular, $c_1(\tilde{\E})\cdot\tilde{H} = 0$. By nefness, $c_1(\tilde{\E})^2 \ge 0$. Hence, by Hodge Index Theorem,
$c_1(\tilde{\E})^2 = 0$
which yields, by \cite[Prop.2.1, Thm.2.5]{DemPetSch}, $c_2(\tilde{\E}) = 0$. So we obtain:
$$c_1(\E)^2 = c_2(\E) = 0.\vspace{-2em}$$
}

\demode{the first assertion in Theorem \ref{main}.}{
Let a variety $X$, an ample $\Q$-Cartier divisor $H$, and a reflexive sheaf $\E$ be as in the asssumptions of Theorem \ref{main}. By Proposition \ref{PseffRestr}, Lemma \ref{Mehta} and \cite[Prop.3.11]{GKPT}, we can consider an integer $m$ big and divisible enough that $n-2$ general members of $|mH|$ cut out a complete intersection normal projective klt surface $S$ in $X$ on which:
\begin{itemize}
\item $\E|_S$ and $(S^{[l]}\E)|_S$, for some $l\ge 6$, are still reflexive;
\item as a consequence, $S^{[l]}(\E|_S)=(S^{[l]}\E)|_S$;
\item both $\E|_S$ and $S^{[l]}(\E|_S)$ remain $H|_S$-stable of zero slope;
\item $\E|_S$ is pseudoeffective.
 \end{itemize}

Then, by Proposition \ref{OnSurface}, there is a finite Galois cover $p:\hat{S}\to S$ such that the reflexive pullback $\hat{\E}:=p^{[*]}\E|_S$ is a nef locally-free sheaf of zero slope. Lemma \ref{NefChern} yields:
$$c_1(\hat{\E})^2 = c_2(\hat{\E}) = 0,$$
so that, by construction, $\hat{c_1}^2(\E|_S)=\hat{c_2}(\E|_S)=0$ and hence:
$$\hat{c_1}^2(\E)\cdot H^{n-2} = \hat{c_2}(\E)\cdot H^{n-2} =0.$$

\noindent The first assertion in Theorem \ref{main} is established.}

\subsection{Proof of Proposition \ref{OnSurface}}\label{subsec:techn}

Let $S$ be a normal projective klt surface, and $H$ an ample $\Q$-Cartier divisor on $S$. Let $\E$ be a reflexive sheaf on $S$ such that:
\begin{itemize}
\item $\hat{c_1}(\E)\cdot H = 0$;
\item $\E$ and $S^{[l]}\E$, for some $l\ge 6$, are stable with respect to $H$;
\item $\E$ is pseudoeffective.
\end{itemize}

As in Section \ref{subsec:orbi}, we denote by $p:\hat{S}\to S$ a finite Galois cover on which the sheaf $\hat{\E}=p^{[*]}\E$ is locally-free. Let $\hat{H}:=p^*H$ be an ample $\Q$-Cartier divisor on $\hat{S}$, $\hat{\pi}:\P(\hat{\E})\to \hat{S}$ be the natural map and $\hat{\zeta}$ be the tautological bundle on $\P(\hat{\E})$.

Abiding by \cite[Sect.3.2]{HoeringPeternell}, we prove two lemmas.
The first lemma uses the stability of $\E$ and of $S^{[l]}\E$, for some $l\ge 6$, to prove the ampleness of $\hat{\zeta}$ on certain subvarieties of $\P(\hat{\E})$.

\lem[Ample]{Keep the notations. For any closed proper subvariety $Z\subset \P(\hat{\E})$ such that the image $\hat{\pi}(Z)$ is not a point in $\hat{S}$, for $m$ big and divisible enough and for a very general curve $\hat{C}\in p^*|mH|$, the restricted tautological $\hat{\zeta}|_{Z\cap \hat{\pi}^{-1}(\hat{C})}$ is ample.}

\demo{Let $Z\subset\P(\hat{\E})$ be a closed proper subvariety whose image $\hat{\pi}(Z)$ has dimension 1 or 2 in $\hat{S}$. Since $p$ is finite, $p(\hat{\pi}(Z))$ has dimension 1 or 2 in $S$. Hence, for $m$ big and divisible enough, a very general curve $C\in |mH|$ satisfies:
\begin{itemize}
\item $C$ is a smooth curve inside the locus $S_0\subset S_{\mathrm{reg}}$ where $\E$ is locally-free;
\item $\hat{C}:=p^{-1}(C)$ is still very general in $p^*|mH|$ and hence smooth too;
\item consequentially, we have locally-free sheaf isomorphisms $\hat{\E}|_{\hat{C}}=p^*\E|_C$ and $S^{l}(\E|_C)=(S^{[l]}\E)|_C$;
\item since $mH$ is ample, $Z\cap \hat{\pi}^{-1}(\hat{C}) \ne \emptyset$; 
\item since $Z$ is proper in $\P(\E)$, $Z\cap \hat{\pi}^{-1}(\hat{C})$ is proper in $\hat{\pi}^{-1}(\hat{C})$;
\item both $\E|_C$ and $S^{6}(\E|_C)$ remain $H|_C$-stable of zero slope, by Lemma \ref{Mehta}.
\end{itemize}

Apply now Lemma \ref{Key} and Remark \ref{QuasiBig}: they establish that $\hat{\zeta}|_{\hat{\pi}^{-1}(\hat{C})}$ is nef and that, for any closed proper variety $W\subset\hat{\pi}^{-1}(\hat{C})=\P(\hat{\E}|_{\hat{C}})$,
$${\left(\hat{\zeta}|_{\hat{\pi}^{-1}(\hat{C})}\right)}^{\,\mathrm{dim}\, W } \cdot W > 0.$$

Using this formula for any closed subvariety $W$ of $Z\cap \hat{\pi}^{-1}(\hat{C})$, the Nakai-Moishezon criterion shows that $\hat{\zeta}|_{Z\cap \hat{\pi}^{-1}(\hat{C})}$ is ample.
}

The second lemma is set at the higher level of $(\hat{S},\hat{\E})$ directly. It uses the pseudoeffectivity and $\hat{H}$-semistability of the locally-free sheaf $\hat{\E}$, infered by Proposition \ref{PseffPullb} and Lemma \ref{SemisPullb}, but no other property of $\E$.

\lem[NotBig]{Keep the notations.
If $\hat{\zeta}$ is not nef, then there is a closed proper subvariety $W$ of $\P(\hat{\E})$ such that,
for $m$ big and divisible enough and for a very general curve $\hat{C}\in p^*|mH|$:
$$\emptyset\ne{W\cap\hat{\pi}^{-1}(\hat{C})}\subsetneq\hat{\pi}^{-1}(\hat{C}),\quad {\hat{\zeta}|_{W\cap\hat{\pi}^{-1}(\hat{C})}}\mbox{ is nef and not big}.$$}

\noindent This result essentially relies on \cite[Lem.3.4]{HoeringPeternell}.

\demo{Denote by $\mu:\tilde{S}\to \hat{S}$ the minimal resolution of $\hat{S}$, by $\tilde{\E}:=\mu^*\hat{\E}$, by $\tilde{\zeta}$ the tautological bundle of $\P(\tilde{\E})$. We have a Cartesian diagram with compatibility of tautological bundles:

\[
\xymatrix{
\P(\tilde{\E}) \ar[d]^*{\tilde{\pi}} \ar[r]^*{\mu'} & \P(\hat{\E})\ar[d]^*{\hat{\pi}}\\
\tilde{S} \ar[r]^*{\mu} & \hat{S}
} 
\]

\noindent Note that $\P(\tilde{\E})$ with its tautological $\tilde{\zeta}$ is a smooth modification of $\P(\hat{\E})$ just as in Definition \ref{PseffNaka}. Hence, $\tilde{\zeta}$ is pseudoeffective.

Suppose that $\hat{\zeta}$ is not nef. In particular, $\tilde{\zeta}$ is not nef, though it is $\tilde{\pi}$-ample. Let $Z\subset B_-(\tilde{\zeta})$ be an irreducible component of maximal dimension. Note that $Z$ contains a $\tilde{\zeta}$-negative curve $N$: its image $\mu'(N)$ must be a $\hat{\zeta}$-negative curve, hence it is not in a fiber of $\hat{\pi}$. So $\hat{\pi}(\mu'(Z))$ is not a point in $\hat{S}$. Moreover, since $\tilde{\zeta}$ is pseudoeffective, $Z\neq\P(\tilde{\E})$.

\medskip

Now, for a very general curve $\hat{C}\in p^*|mH|$ for $m$ big and divisible enough, 
\begin{itemize}
\item $\hat{C}\subset \hat{S}_{\mathrm{reg}}$; in particular, $\mu$ is an isomorphism over $\hat{C}$;
\item $\emptyset\ne Z\cap\mu'^{-1}(\hat{\pi}^{-1}(\hat{C}))\subsetneq \mu'^{-1}(\hat{\pi}^{-1}(\hat{C}))$; 
\item $\hat{\E}|_{\hat{C}}$ is nef by Lem.\ref{Mehta} and Prop.\ref{Miya} and it has $\hat{H}|_{\hat{C}}$-slope zero;
\item hence, $\tilde{\zeta}|_{\mu'^{-1}(\hat{\pi}^{-1}(\hat{C}))}$ is nef too, and moreover its top power is zero;
\item hence, by \cite[Lem.3.4]{HoeringPeternell} (which applies since $Z$ was chosen with minimal codimension):
	\begin{align*}0 &= {\left(\tilde{\zeta}|_{\mu'^{-1}(\hat{\pi}^{-1}(\hat{C}))}\right)}^{\,\mathrm{dim}\,\mu'^{-1}(\hat{\pi}^{-1}(\hat{C}))}\\
		&\ge {\left(\tilde{\zeta}|_{Z\cap\mu'^{-1}(\hat{\pi}^{-1}(\hat{C}))}\right)}^{\,\mathrm{dim}\,Z\cap\mu'^{-1}(\hat{\pi}^{-1}(\hat{C}))}
\ge 0.\end{align*}
\end{itemize}

As $\mu'$ is an isomorphism over $\hat{C}$, $W:=\mu'(Z)$ works well as the closed proper subvariety of $\P(\hat{\E})$ we wanted to construct.
}

We now combine these lemmas to establish Proposition \ref{OnSurface}.

\demode{Proposition \ref{OnSurface}.}{
Suppose by contradiction that $\hat{\E}$ is not nef. Then Lemma \ref{NotBig} yields a closed proper subvariety $W$ of $\P(\hat{\E})$ which satisfies, for $m$ big and divisible enough and for a very general curve $\hat{C}\in p^*|mH|$:

$$\emptyset\ne W\cap \hat{\pi}^{-1}(\hat{C})\subsetneq \hat{\pi}^{-1}(\hat{C})\mbox{ and $\hat{\zeta}|_{W\cap \hat{\pi}^{-1}(\hat{C})}$ is nef and not big}.$$
The first condition shows that $\hat{\pi}(W)$ is not a point. So Lemma \ref{Ample} applies, hence $\hat{\zeta}|_{W\cap \hat{\pi}^{-1}(\hat{C})}$ is ample, contradiction!
}

\section{Proof of Theorem \ref{main}}

As it follows from the discussion in Section \ref{subsec:surf}, Theorem \ref{main} is halfway. Here is what remains to prove:

\theo[LuTa]{Let $X$ be a normal projective klt variety of dimension $n$ with an ample $\Q$-Cartier divisor $H$. Let $\E$ be a reflexive sheaf on $X$, such that:
\begin{itemize}
\item $\E$ is $H$-semistable;
\item the following equalities hold:
$$\hat{c_1}(\E)\cdot H^{n-1} = \hat{c_1}^2(\E)\cdot H^{n-2}
= \hat{c_2}(\E)\cdot H^{n-2} = 0.$$
\end{itemize}
\noindent Then there is a finite Galois morphism $\nu:\tilde{X}\to X$, étale in codimension 1, such that $\nu^{[*]}\E$ is a locally-free sheaf with numerically trivial determinant, and is $\mathrm{Gal}(\tilde{X}/X)$-equivariantly flat. 
Consequentially, $\nu^{[*]}\E$ is numerically flat and its first and second Chern classes are numerically trivial.}

\demo{We apply \cite[Thm.1.4]{LuTaji} to obtain a finite Galois morphism $\nu:\tilde{X}\to X$, étale over $X_{\mathrm{reg}}$, such that $\nu^{[*]}\E$ is locally-free with a numerically trivial determinant and $\mathrm{Gal}(\tilde{X}/X)$-equivariantly flat. 

Let then $\eps:\tilde{X}'\to\tilde{X}$ be a resolution of $\tilde{X}$ and $\E':=\eps^*\nu^{[*]}\E$, which is a flat locally-free sheaf with a numerically trivial determinant on $\tilde{X}'$.
As shown in \cite[Rmk.2.6]{HoeringPeternell}, $\E'$ is then numerically flat and its Chern classes vanish (as cohomological classes on $\tilde{X}'$). By Prop.\ref{moreonnef}, $\nu^{[*]}\E$ is nef, hence numerically flat.
Moreover, for any $\Q$-Cartier divisors $D_1,\ldots, D_{n-2}$,
$$c_2(\nu^{[*]}\E)\cdot D_1\cdots D_{n-2}
= c_2(\E')\cdot \eps^*D_1\cdots \eps^*D_{n-2}
= 0,$$
\noindent so the Chern classes of $\nu^{[*]}\E$ are trivial, which completes the proof of the theorem.}

\section{Proof of Theorem \ref{CYNotPseff}}\label{sec:thmnotpseff}

We give a few definitions along the lines of Theorem \ref{BBDec}:

\defi[CY]{Let $X$ be a normal projective canonical variety of dimension $n\ge 2$. It is called:
\begin{itemize}
\item a {\it Calabi-Yau} variety if $h^0(Y,\Omega^{[q]}_Y)=0$ for all integers $1\le q\le n-1$ and all quasiétale finite covers $Y\to X$;
\item an {\it irreducible holomorphic symplectic (IHS)} variety if there is a reflexive form $\sigma\in H^0(X,\Omega^{[2]}_X)$ such that, for any quasiétale finite cover $f:Y\to X$, the reflexive form $f^{[*]}\sigma$ generates $H^0(X,\Omega^{[\cdot]}_Y)$ as an algebra for the wedge product.
\end{itemize}
}

We use the terms singular Calabi-Yau (resp. IHS) variety and Calabi-Yau (resp. IHS) variety interchangeably, unless explicitly said otherwise. They may both accidentally denote smooth varieties.

\defi[K3]{For the sake of a consistent terminology, let us call a {\it singular K3 surface}, or for short a {\it K3 surface}, a normal projective klt surface which has no finite quasiétale cover by an abelian variety. Equivalently, it is a Calabi-Yau variety or an IHS variety of dimension 2.}

\defi{For the sake of a convenient vocabulary, let us define the {\it augmented irregularity} $\tilde{q}(X)$ of a normal projective klt variety $X$ with trivial canonical class as the maximum of all irregularities $q(Y)$ of finite quasiétale covers $Y$ of $X$. Note that it is precisely the dimension of the abelian part in the singular Beauville-Bogomolov decomposition of $X$.}

\medskip

Let us now proceed to prove Theorem \ref{CYNotPseff}.

\demode{Theorem \ref{CYNotPseff}. }{Let $X$ be a normal projective klt variety of dimension at least 2 with trivial canonical class. Suppose that $\Omega^{[1]}_X$ is pseudoeffective (the same whole argument works just alike for the tangent sheaf $\mathcal{T}_X$) and assume by contradiction that $\tilde{q}(X)=0$.

The singular Beauville-Bogomolov decomposition then reads:
$$f:\tilde{X}\to X\mbox{ and }\tilde{X}\cong\prod_i Y_i\times\prod_j Z_j,$$
with the same notations as in Theorem \ref{BBDec}.

Remember that $f^{[*]}\Omega^{[1]}_X = \Omega^{[1]}_{\tilde{X}}$, since reflexive sheaves are normal and there is a big open set over which $f$ is just a finite étale cover. By Proposition \ref{PseffPullb}, $\Omega^{[1]}_{\tilde{X}}$ is pseudoeffective; it splits according to the product defining $\tilde{X}$. So there is a factor $Y$ (Calabi-Yau or IHS) of $\tilde{X}$ such that $\Omega^{[1]}_{Y}$ is pseudoeffective \cite[inductive argument in Proof of Thm.1.6]{HoeringPeternell}. Now, $\Omega^{[1]}_{Y}$ satisfies all hypotheses of Theorem \ref{main}, the stability assumptions coming from \cite[Prop.8.20]{GKPSing} and \cite[Rmk.8.3]{GGK}. 

As a consequence, for some ample polarization $H$ on $Y$, $\hat{c_2}(\Omega^{[1]}_Y)\!\cdot\! H^{\mathrm{dim}\, Y}\!=\!0$, so that $Y$ has a finite quasiétale cover by an abelian variety by \cite[Thm.1.4]{LuTaji}, contradiction!}

\rem{This pseudoeffectiveness result can be considered as an interesting improvement of the effectiveness result \cite[Thm.11.1]{GGK}, which says that $\tilde{q}(X)=0$ if and only if, for all $m\in\N$, $h^0(X,S^{[m]}\Omega^{[1]}_X) = 0$.}

\medskip

Examples for Theorem \ref{CYNotPseff} are to search among normal projective klt varieties with trivial canonical class singularities in codimension 2, which are plethoric. But singular varieties whose decomposition is known are not so numerous; and, for sure, one shall understand the Beauville-Bogomolov type of a given variety before telling anything about the positivity of its reflexivized cotangent sheaf.

\expl[GlobalQ]{A first example to which Theorem \ref{CYNotPseff} applies is the following \cite[Par.14.2.2]{GGK}: let $F$ be a Fano manifold on which a finite group $G$ acts freely in codimension 1. Suppose there is a smooth $G$-invariant element $Y$ in the linear system $|-K_F|$. Then, $Y$ is a smooth Calabi-Yau variety with a $G$-action. If the volume form on $Y$ is preserved by this action, then $X:=Y/G$ is a normal projective klt variety with trivial canonical class, and the morphism $Y\to X$ has no ramification divisor, hence it is quasiétale. The fact that the decomposition of $X$ consists of a smooth Calabi-Yau manifold $Y$ guarantees that $X$ is a singular Calabi-Yau variety, as presented in Definition \ref{CY}.

Although $X$ may well have singularities in codimension 2, they merely stem from its global quasiétale quotient structure. In particular, \cite[Thm.1.6]{HoeringPeternell} actually proves the non-pseudoeffectiveness of $\mathcal{T}_X$ and $\Omega^{[1]}_X$, namely because it applies to $Y$ and converts onto $X$ through Proposition \ref{PseffPullb}. Hence, the example is quite shallow: it has no real need for the machinery dealing with singularities in codimension 2 that Theorem \ref{CYNotPseff} is about.}

In the next section, we present better examples for Theorem \ref{CYNotPseff}, namely Calabi-Yau threefolds with singularities in codimension 2 that are not constructed as global quasiétale quotients of varieties which are smooth in codimension 2.

\section{Threefolds in Theorems \ref{BBDec} and \ref{CYNotPseff}}\label{sec:threeexpl}

In Section \ref{sec:thmnotpseff}, we defined singular Calabi-Yau and IHS varieties. They can also be defined by means of their algebraic holonomy, an approach which \cite{GGK} uses thoroughly. This notably enables one to prove that IHS varieties must have even dimension \cite[Thm.12.1, Prop.12.10]{GGK}. In particular, the singular Beauville-Bogomolov decomposition for a normal projective klt variety $X$ of dimension 3 is quite simple:
$\tilde{X}$ has to be one of the following:
\begin{itemize}
\item a smooth abelian variety;
\item a product $S\times E$, where $S$ is a K3 surface as in Definition \ref{K3} and $E$ is a smooth elliptic curve;
\item a Calabi-Yau variety.
\end{itemize}

The aforementioned \cite[Thm.1.4]{LuTaji} provides a criterion for identifying the purely abelian case by computing $\hat{c_2}(X)$. 

One is then left with two cases: the singular threefold $X$ may arise from a product $S\times E$, in which case $\mathcal{T}_X$ and $\Omega^{[1]}_X$ are pseudoeffective because of the abelian factor $E$; alternatively, $X$ can be a genuine singular Calabi-Yau threefold. This second possibility is hard to identify, but, when it happens, it may give new examples for Theorem \ref{CYNotPseff}. 

The next subsection is devoted to providing a necessary condition for a normal projective klt threefold to be finitely quasiétaly covered by a product $S\times E$.

\subsection{Products of a K3 surface and an elliptic curve}

We are going to prove the following result:

\prop[pasK3quotient]{Let $X$ be a normal projective klt threefold with trivial canonical class. Suppose its Beauville-Bogomolov decomposition is of the form 
$$\tilde{X} = S\times E,$$
where $S$ is a K3 surface and $E$ a smooth elliptic curve.
Then $X$ has fibrations:

\[\xymatrix{
 X \ar[d] \ar[r] & S/G_S\\
E/G_E
} 
\]
where $G_E$ and $G_S$ are finite subgroups of $\mathrm{Aut}(E)$ and $\mathrm{Aut}(S)$.
 In particular, $\rho(X)\ge 2$.}

Let us first state a weak uniqueness result, guaranteeing that the statement of Proposition \ref{pasK3quotient} makes sense. It is straightforward from the proof of the Beauville-Bogomolov decomposition theorem.

\prop[unicity]{Let $X$ be a normal projective klt variety with trivial canonical class. Then the number, types and dimensions of the factors of a finite quasiétale covering $\tilde{X}\to X$ as in Theorem \ref{BBDec} do not depend on the choice of that covering.}

A finite quasiétale morphism is not necessarily a quotient map by a finite group action free in codimension 1. In the smooth case however, \cite[Lem.p.9]{Beauv} allows us to assume that the finite étale decomposition morphism $p:\tilde{X}\to X$ is Galois. Let us state a partial singular analog:

\prop[Galois]{Let $X$ be a normal projective klt variety with trivial canonical class. Take a finite quasiétale covering $f:\tilde{X}\to X$ as in Theorem \ref{BBDec}.
Suppose that all Calabi-Yau factors of $\tilde{X}$ have even dimension. Then there is a finite quasiétale Galois morphism $f':\tilde{Z}\to X$, so that $\tilde{Z}$ splits into factors in the same number, types, and dimensions as $\tilde{X}$.} 

\demo{By \cite[Thm.1.5]{GKPFund}, we can take a finite quasiétale Galois covering $g:Y\to X$ such that any finite morphism $Z\to Y$ étale over $Y_{\mathrm{reg}}$ is étale over $Y$. By purity of the branch locus, any quasiétale morphism $Z\to Y$ is then étale.

Note that $Y$ is still a normal projective klt variety with trivial canonical class, hence has a singular Beauville-Bogomolov decomposition $h:Z\to Y$. By Proposition \ref{unicity}, the factors of $Z$ have the same type as those of $\tilde{X}$. It writes:

$$Z=A\times\prod_i Y_i\times\prod_j Z_j,$$ 
where $A$ is an abelian variety, $Y_i$ Calabi-Yau varieties and $Z_j$ IHS varieties. Since all $Y_i$ and, of course, all $Z_j$ have even dimension, by \cite[Cor.13.3]{GGK}, they are simply connected. 

Hence, finite étale fundamental groups equal: $\widehat{\pi_1}(Z)\simeq\widehat{\pi_1}(A)$. That is to say, any finite étale cover of $Z$ actually stems from a finite étale cover of $A$.

We now use \cite[Thm.3.16]{GKPFund}: there is a finite Galois morphism $\gamma:\tilde{Z}\to Z$ such that $\Gamma=g\circ h\circ \gamma:\tilde{Z}\to X$ is finite Galois and ramifies where $g\circ h$ does. So $\Gamma$ is still quasiétale, in particular $h\circ\gamma:\tilde{Z}\to Y$ is quasiétale too. By construction of $Y$, $h\circ\gamma$ is then étale, so that $\gamma$ is étale. By construction of $Z$, one has:

$$\tilde{Z}=A'\times\prod_i Y_i\times\prod_j Z_j,$$
where $A'$ is a finite étale cover of the abelian variety $A$.
Finally, $\Gamma:\tilde{Z}\to X$ is finite Galois quasiétale, and $\tilde{Z}$ splits as mandated.}

\rem{The main obstacle for generalizing this proposition is the fact that fundamental groups of odd-dimensional Calabi-Yau varieties are poorly understood \cite[Sect.13.2]{GGK}; most notably, they may not be finite.
}

Here is the last ingredient for the proof of Proposition \ref{pasK3quotient}:

\lem[auto]{Let $S$ be a K3 surface as in Definition \ref{K3}, $E$ a smooth elliptic curve. Then:
$$\mathrm{Aut}(S\times E)\cong\mathrm{Aut}(S)\times\mathrm{Aut}(E).$$}

\demo{Let $\tilde{S}$ be the minimal resolution of $S$. It is a smooth K3 surface, so $\mathrm{Aut}(\tilde{S})$ is discrete. Moreover, the uniqueness of minimal resolution implies that any automorphism of $S$ lifts to an automorphism of $\tilde{S}$, and this is obviously an injection. Hence, $\mathrm{Aut}(S)$ is discrete.

Let us now copy the argument by \cite[Lem.p.8]{Beauv}. Let $u\in\mathrm{Aut}(S\times E)$. Since the projection $p_E:S\times E\to E$ is the Albanese map of $S\times E$, we can factor $p_E\circ u$ by it: there is $v\in\mathrm{Aut}(E)$ such that $p_E\circ u = v\circ p_E$.
Hence, there is a map $w:E\to\mathrm{Aut}(S)$ which decomposes:
$$u:(s,e)\in S\times E\mapsto (w_e(s),v(e)).$$
Since $\mathrm{Aut}(S)$ is discrete, the map $w$ is constant, so $u=(w_0,v)$.}

\demode{Proposition \ref{pasK3quotient}. }{Let $X$ be a normal projective variety of dimension 3 with trivial canonical class. Suppose that there is a finite quasiétale cover $f:S\times E\to X$, where $S$ is a singular K3 surface and $E$ a smooth elliptic curve. By Proposition \ref{Galois}, we can assume that there is a finite group $G$ acting on $S\times E$ such that $f$ is the induced quotient map. By Lemma \ref{auto}, $G$ can be considered a subgroup of $\mathrm{Aut}(S)\times\mathrm{Aut}(E)$. As it acts diagonally, we have the following diagram:

\[
\xymatrix{
S\times E \ar[r]^*{p_S} \ar[d]^*{p_E} \ar[rd]^*{f} & S \ar[rd]\\
E \ar[rd] & X \ar[d] \ar[r] & S/G_S\\
& E/G_E
} 
\]

\noindent so that $\rho(X)$ is at least 2.}

\subsection{Calabi-Yau hypersurfaces in weighted projective spaces}\label{sec:expl}

The aim of this last part is to provide examples of Calabi-Yau threefolds that are singular along curves, by establishing the following result:

\prop[areCY]{Let $\P=\P(w_0,\ldots,w_{4})$ be a weighted projective space and $d = w_0 + \ldots + w_{4}$ such that there is a general wellformed quasismooth hypersurface $X$ of degree $d$ in $\P$. Suppose that $X$ contains no edge of $\P$. Then $X$ is a singular Calabi-Yau in the sense of Definition \ref{CY}.}

A general exposition to complete intersections in weighted projective spaces can be found in \cite{Fletcher}. We stick to its terminology. 

Let $\P=\P(w_0,\ldots,w_{4})$ be a wellformed 4-dimensional weighted projective space. There is a ramified quotient map: $p:\P^n\to \P$, by the finite diagonal group action of $\bigoplus_i\Z_{w_i}$ on $\P^n$. With homogeneous coordinates on either side, we can write: $$p:[x_0:\ldots:x_n]\in\P^n\mapsto[y_0=x_0^{w_0}:\ldots:y_n=x_n^{w_n}]\in\P.$$

We denote by $\mathcal{O}_{\P}(1)$ the ample $\Q$-Cartier divisor on $\P$ whose pullback by $p$ is $\mathcal{O}_{\P^n}(1)$. If the linear system $|\mathcal{O}_{\P}(d)|$ contains a wellformed quasismooth hypersurface, it actually contains a Zariski-open set of such hypersurfaces and we write $X_d$ for a general one.

\subsubsection*{Singularities of general quasismooth hypersurfaces of dimension 3}
Let $X$ be a general quasismooth hypersurface of degree $d$ and of dimension 3 in the weighted projective space $\P$. Then $X$ is a full suborbifold of $\P$ (see \cite[Def.5]{Borz} for a definition, \cite[Thm.3.1.6]{Dolga} for a proof). In particular, $X_{\mathrm{sing}} = X\cap\P_{\mathrm{sing}}$, and at any point $x\in X\cap\P_{\mathrm{sing}}$, writing that $\P$ is locally isomorphic to a quotient $\C^4/G_x$, $X$ is locally isomorphic to $\C^3/G_x$ in a compatible way with inclusions.
Hence, $X$ has only quotient singularities, so it is klt. The locus $X_{\mathrm{sing}}$ is a finite union of curves and points, which may be of various types:
\begin{itemize}
\item a {\it vertex} in $\P$ is a point with $y_i=1$ for a single $i\in\lint 0,4\rint$ and $y_j=0$ for all $j\ne i$. If $w_i\ne 1$, this vertex is a singular point in $\P$. It gives rise to a singular point in $X$ if and only if it lies in it, ie $w_i$ does not divide $d$.
\item an {\it edge} in $\P$ is a line with equation $y_j = 0$ for all $j\in J$, for a certain $J\subset\lint 0,4\rint$ of cardinal 3. If $\gcd(w_j)_{j\not\in J}\ne 1$, the edge is in $\P_{\mathrm{sing}}$. Recall that $X$ is taken general in its linear system. Hence, an edge in $\P$ lies entirely in $X$ if and only if $(w_j)_{j\not\in J}$ do not partition $d$, in $X_{\mathrm{sing}}$ if and only if $(w_j)_{j\not\in J}$ do not partition $d$ and have a non-trivial common divisor. If an edge in $\P_{\mathrm{sing}}$ does not lie entirely in $X$, it gives a finite amount of points in $X_{\mathrm{sing}}$.
\item a {\it 2-face} in $\P$ is a 2-plane with equation $y_j = 0$ for all $j\in J$, for a certain $J\subset\lint 0,4\rint$ of cardinal 2. If $\gcd(w_j)_{j\not\in J}\ne 1$, the 2-face is in $\P_{\mathrm{sing}}$. By quasismoothness, no 2-face lies entirely in $X$. Hence, any 2-face intersects $X$ along an effective 1-cycle. In this way, 2-faces in $\P_{\mathrm{sing}}$ may produce curves in $X_{\mathrm{sing}}$.
\end{itemize}

Under the additional hypothesis that $X$ contains no edge of $\P$, we can say more about singular loci. 

Indeed, it is worth noticing that the restricted quotient map $p^{-1}(X)\to X$ is an unfolding of $X$, as defined in Section \ref{subsec:orbi}; we may write $\hat{X}$ for $p^{-1}(X)$. For establishing Prop.\ref{areCY}, we will prove the following:

\lem[c2wps]{Let $X$ be a general wellformed quasismooth hypersurface of dimension 3 in a weighted projective space $\P$ not isomorphic to $\P^{4}$. Assume that $X$ has trivial canonical class and that it contains no edge of $\P$. Then $\hat{c_2}(X)\cdot\mathcal{O}_X(1)>0$.}

\noindent In the course of the proof of this lemma, we will use the fact that $X$ containing no edge of $\P$, $\hat{X}$ is smooth in codimension 2.

\rem{Note that the restricted finite map $\hat{X}=p^{-1}(X)\to X$ is certainly ramified along divisors, so that $X$, in the lucky case where it happens to be a singular Calabi-Yau threefold, is not at all constructed as a finite quasiétale global quotient, contrarily to the unsatisfying Example \ref{GlobalQ}.}

The proof that $\hat{X}$ is smooth in codimension 2 relies on the following lemma and remark:

\lem[base]{Let $X$ be a general quasismooth hypersurface of degree $d$ in the weighted projective space $\P=\P(w_0,\ldots,w_4)$. Suppose that it contains no edge of $\P$. Then the base locus $\mathrm{Bs}(\mathcal{O}_{\P}(d))$ has dimension 0.}

\demo{Let $Z$ be an irreducible component of the base locus of $\mathcal{O}_{\P}(d)$, let us prove by induction on $\dim\P$ that it is a point. Suppose we are at the induction step where the ambient space $\P'$ has local coordinates $y_0,y_1,y_2,\ldots$ and dimension 4, 3 or 2.

Denote by $H_i$ the hyperplane $\{y_i=0\}$ in $\P'$, by $\P'_i$ the isomorphic\linebreak weighted projective space $\P'(\ldots,\hat{w_i},\ldots)$. By \cite[Prop.4.A.3]{Beltr}, we have an isomorphism between the restriction $\mathcal{O}_{\P'}(d)\otimes\mathcal{O}_{H_i}$ and the $\Q$-Cartier divisor $\mathcal{O}_{\P'_i}(d)$. This translates to global sections as a surjection:
\begin{equation}\label{sections}
H^0(\P',\mathcal{O}_{\P'}(d))\twoheadrightarrow H^0(\P'_i,\mathcal{O}_{\P'_i}(d)),
\end{equation}
which is given by setting $y_i=0$ when considering the global sections as certain polonomials in the local coordinates of $\P'$.

The quasismoothness of $X$ in $\P$ and the way the composite surjection 
$$H^0(\P,\mathcal{O}_{\P}(d))\twoheadrightarrow H^0(\P',\mathcal{O}_{\P'}(d)),$$ 
writes in local coordinates yield a global section of $\mathcal{O}_{\P'}(d)$ of the form $y_0^{\alpha_0}y_1^{\alpha_1}y_2^{\alpha_2}$. In particular, there is an $i=0,1$ or $2$ such that $Z\subset H_i\simeq\P'_i$.
Moreover, by Eq.\ref{sections}, $Z$ sits in the base locus of $\mathcal{O}_{\P'_i}(d)$. 

Induction propagates from $\P'=\P$ down to when we obtain that $Z$ is contained in an edge $H_{ijk}$ of $\P$ and in the base locus $\mathrm{Bs}(\mathcal{O}_{\P_{ijk}}(d))\subset\mathrm{Bs}(\mathcal{O}_{\P}(d))\subset X$. Since $X$ contains no edge of $\P$, $Z$ is in $X\cap H_{ijk}$ of dimension 0, so it is a point.}

\rem[singinter]{With the same notations and hypotheses, the intersection of $X$ with any 2-face of $\P_{\mathrm{sing}}$ is a reduced curve.}

\demo{As in the proof of Lemma \ref{base}, the intersection is scheme-theoretically defined by a general section of $\mathcal{O}_{\P_{ij}}(d)$. We are to show that such general section of $\mathcal{O}_{\P_{ij}}(d)$ is quasismooth in the weighted projective space $\P_{ij}$, hence it is a variety by \cite[3.1.6]{Fletcher}.

We use the arithmetical criterion for quasismoothness: since $X$ contains no edge of $\P$, each pair $w_a,w_b$ partitions $d$. We are left to check the criterion for $k=1$: fix any $a\ne i,j$, we want to find $b\ne i,j$ such that $w_a$ divides $d-w_b$. It is clear that there is a $b\in\lint 0,4\rint$ satisfying that. As $H_{ij}$ is a 2-face in $\P_{\mathrm{sing}}$, the greatest common divisor of all weights except $w_i,w_j$ is non-trivial, divides $d$ but neither $w_i$ nor $w_j$ (by wellformedness). In particular, since this greatest common divisor divides $w_b=d-\alpha w_a$, $b\ne i,j$, as wished.}

We can now deduce:

\prop[unfoldsing]{Let $X$ be a general quasismooth hypersurface of degree $d$ in a weighted projective space $\P=\P(w_0,\ldots,w_4)$, $p$ the natural quotient $\P^4\to\P$, $\hat{X}=p^{-1}(X)$. Suppose that $X$ contains no edge of $\P$. Then $\hat{X}$ is smooth in codimension 2.}

\demo{The threefold $\hat{X}$ is general in the linear system $p^*|\mathcal{O}_{\P}(d)|$, whose base locus has dimension 0 by Lemma \ref{base}. By Bertini's theorem, $\hat{X}$ is smooth in codimension 2.
}

\rem{The converse of Proposition \ref{unfoldsing} does not hold: for instance, the general quasismooth $X_7$ in $\P(1,1,1,2,2)$ contains the edge of equation $y_0=y_1=y_2=0$, but its unfolding is nevertheless smooth in codimension 2.}

\expl{The hypothesis of Proposition \ref{unfoldsing} is not that $X$ contains no edge of $\P_{\mathrm{sing}}$, but that it contains no edge of $\P$ at all: for instance, consider the general $X=X_{56}$ in $\P(2, 4, 9, 13, 28)$. It contains a single edge of $\P$, namely $e$ of equation $y_0=y_1=y_4=0$. This edge does not actually lie in $\P_{\mathrm{sing}}$, as 9 and 13 are coprime, but one can check that $\hat{X}$ has the curve $p^{-1}(e)$ in its singular locus (by computing the derivatives of the equation defining $\hat{X}$ in $\P^4$ along the curve $p^{-1}(e)$).}

\expl[ex1734]{The general wellformed quasismooth hypersurface $X=X_{1734}$ in $\P(91,96,102,578,867)$ contains no edge of $\P$. In particular, $\hat{X}$ is smooth in codimension 2 by Proposition \ref{unfoldsing}. 

Moreover, the curves of $X_{\mathrm{sing}}$ are precisely the intersections of $X$ with all 2-faces of $\P_{\mathrm{sing}}$, which we can list:
\begin{itemize}
\item $y_0=y_1=0$ of type $\frac{1}{17}(6,11)$,
\item $y_0=y_3=0$ of type $\frac{1}{3}(1,2)$,
\item $y_0=y_4=0$ of type $\frac{1}{2}(1,1)$.
\end{itemize}
}

It is possible to check the type of singularities of a general hypersurface of a given degree in a given weighted projective space by a simple computer program.

\subsubsection*{Proof of Proposition \ref{areCY}} As we said before, the main ingredient in the proof is Lemma \ref{c2wps}.

\demode{Lemma \ref{c2wps}.}{Let $p:\P^{4}\to \P$ be the natural quotient map. Writing $\P = \P(w_0,\ldots,w_4)$ with $(w_0,\ldots,w_4)$ not colinear to $(1,\ldots,1)$, the morphism $p$ has degree $w_0\cdots w_4$, which we denote by $N$, and $X$ has degree $w_0+\ldots+w_4$, which we denote by $d$. We may also write $s$ for the symmmetric elementary polynomial of degree 2 in the weights and $q$ for the sum of their squares: $d^2=q+2s$.

Since $X$ is a full suborbifold of $\P$, $\hat{X}:=p^{-1}(X)\to X$ is an unfolding of $X$ as defined in Section \ref{subsec:orbi}. Applying the left-exact functor of reflexive pullback (see Lemma \ref{ReflPB}) to the exact sequence:
$$0\to\mathcal{T}_X\to\mathcal{T}_{\P}|_X \to -K_{\P},$$
we get another exact sequence:
$$0\to p^{[*]}\mathcal{T}_X \to p^{[*]}\mathcal{T}_{\P}|_X \to p^{[*]}(-K_{\P}) \to \mathcal{Z}\to 0,$$
where the coherent sheaf $\mathcal{Z}$ is supported on the locus $p^{-1}(\mathrm{Sing}\,X)\subset\hat{X}$ of codimension at least 2.

Because of the last surjection, $\dim_{k(p)}\mathcal{Z}\otimes \mathcal{O}_p\le 1$ for any closed point $p\in \tilde{X}$. 

By Proposition \ref{unfoldsing}, the unfolding $\hat{X}$ is smooth in codimension 2, so the usual second Chern class $c_2(\mathcal{Z})$ makes sense. Since usual Chern classes are additive, and $c_1(\mathcal{T}_{X})=0,c_1(\mathcal{Z})=0$:
\begin{equation*}
\hat{c_2}(\mathcal{T}_X)\cdot\mathcal{O}_X(1) = 
\hat{c_2}(\mathcal{T}_{\P}|_X)\cdot\mathcal{O}_X(1) 
+\frac{1}{N} c_2(\mathcal{Z})\cdot\mathcal{O}_{\hat{X}}(1).
\end{equation*}

By the Miyaoka-Yau inequality \cite[Thm.1.5]{GKPPreprint}, we have a positive contribution:
$$\hat{c_2}(\mathcal{T}_{\P}|_X)\cdot \mathcal{O}_X(1)
= \hat{c_2}(\mathcal{T}_{\P})\cdot (-K_{\P})\cdot\mathcal{O}_{\P}(1)
\ge \frac{4}{10}(-K_{\P})^3\cdot\mathcal{O}_{\P}(1)
= \frac{4d^3}{10N}.$$

Let us estimate the other summand. Take $m$ big and divisible enough that 
$\mathcal{O}_{\hat{X}}(m)$ is very ample and $S$ a general element in $|\mathcal{O}_{\hat{X}}(m)|$. By \cite[Lem.10.9]{FlipsNA},
$$c_2(\mathcal{Z})\cdot \mathcal{O}_{\hat{X}}(1) 
= \frac{1}{m}c_2(\mathcal{Z}|_S)
= -\frac{1}{m}\mathrm{deg}(\mathcal{Z}|_S)$$

Denote by $C_1,\ldots C_k$ the curves in $X_{\mathrm{sing}}$. By Lemma \ref{sublemma}, we can bound:
\begin{align*}
\mathrm{deg}(\mathcal{Z}|_S)
\le\mathrm{Card}\left(S\cap\bigcup_{i=1}^k p^{-1}(C_i)\right)
&=\sum_{i=1}^k N\mathcal{O}_X(m)\cdot C_i\\
&\le Nm\mathcal{O}_X(1)^3\sum_{0\le i<j\le 4} w_iw_j\\
&= mNs (-K_{\P})\cdot\mathcal{O}_{\P}(1)^3\\
&= msd.
\end{align*}

Finally putting the positive and negative part together,
\begin{align*}
\hat{c_2}(X)\cdot\mathcal{O}_X(1)
&> \frac{4d^3-10sd}{10N}\\
&= \frac{d(4q-2s)}{10N}\\
&=\frac{d}{10N}\sum_{0\le i<j\le 4}(w_i-w_j)^2
> 0.
\end{align*}}

\lem[sublemma]{Let $X$ be a general wellformed quasismooth hypersurface of dimension 3 in a weighted projective space $\P$. Assume that $X$ has trivial canonical class and contains no edge of $\P$. Then there are at most 10 curves in $X_{\mathrm{sing}}$, with different cohomological classes in the list of the $$[\mathcal{O}_X(w_i)\cdot\mathcal{O}_X(w_j)]\in H^4(X;\Q),\mbox{ for }0\le i< j\le 4.$$}

\demode{Lemma \ref{sublemma}.}{By Remark \ref{singinter}, each curve in $X_{\mathrm{sing}}$ is scheme-\linebreak theoretically the complete intersection of $X$ with a 2-face $H_{ij}$ of $\P_{\mathrm{sing}}$. This association being bijective, there are as many curves in $X_{\mathrm{sing}}$ as 2-faces in $\P_{\mathrm{sing}}$, so at most 10. The curve that corresponds to the 2-face $H_{ij}$ has cohomological class $[\mathcal{O}_X(w_i)\cdot\mathcal{O}_X(w_j)]$.}

Now we can finally establish Proposition \ref{areCY}:

\demo{Consider $X$ a general wellformed quasismooth hypersurface of degree $d=w_0+\ldots+w_4$ in a weighted projective space $\P=\P(w_0,\ldots, w_4)$. Suppose that $X$ contains no edge of $\P$.
If $\P$ is $\P^4$, $X$ is smooth and there is nothing to prove. Let us assume $\P\not\cong\P^4$. By Lemma \ref{c2wps}, $\hat{c_2}(X)\cdot\mathcal{O}_X(1)\ne 0$, hence by \cite[Thm.1.4]{LuTaji}, $X$ is not a finite quotient of an abelian threefold. Moreover, one has $\pic(X)\simeq\Z$ \cite[Thm.3.2.4(i)]{Dolga}, so Proposition \ref{pasK3quotient} applies to $X$: it is not covered by a product of a K3 surface and an elliptic curve, hence its Beauville-Bogomolov decomposition consists of a single Calabi-Yau factor. By Lemmma \ref{fundgroup}, so $X$ itself is a Calabi-Yau variety, in the sense of Definition \ref{CY}. In particular, $X$ has canonical (and not merely klt) singularities.}

\lem[fundgroup]{Let $X$ be a general quasismooth hypersurface in a weighted projective space $\P$. Then any finite quasiétale cover $X'$ of $X$ is trivial.}

\demo{Let $X'$ be a finite quasiétale cover of $X$ of degree $d$; note that by Zariski purity of branch locus, it is étale over $X_\mathrm{reg}$. Let $C_X^*\subset\C^{n+1}\setminus\{0\}$ be the smooth cone over $X$, with the projection $q:C_X^*\to X$. 
The morphism $C'=X'\underset{X}{\times}C_X^*\to C_X^*$ is finite of degree $d$ and étale over the big open set $q^{-1}(X_{\mathrm{reg}})\subset C_X^*$. 
Normalizing, the map $\tilde{C'}\to C_X^*$ has degree $d$ and is étale over a big open set as well. 
As $C_X^*$ is smooth, this map is actually finite étale; by \cite[Lem.3.2.2(ii)]{Dolga}, 
$\pi_1^{\mbox{{\footnotesize ét}}}(C_X^*)=\{1\}$ so $d=1$.
}

\subsubsection*{Examples for Proposition \ref{areCY}} General wellformed quasismooth hypersurfaces with trivial canonical class in 4-dimensional weighted projective spaces are classified in \cite{KreuzerSkarke}. There is an explicit exhaustive list of the 7555 of them. In this list, 7238 elements are not smooth in codimension 2, and 2409 elements that are not smooth in codimension 2 also contain no edge of their ambient weighted projective space. These elements fulfill the hypotheses for Proposition \ref{areCY}, just as Example \ref{ex1734} did: so they are singular Calabi-Yau threefolds to which Theorem \ref{CYNotPseff} applies.

The exhaustive enumerations of elements of the \cite{KreuzerSkarke} classification satisfying additional properties were done by running a simple computer program on the database \cite{data}.

\rem{For the reader misguided by variations in terminology, the arXiv version of this paper proves the elementary fact that the varieties studied in \cite{KreuzerSkarke} are the same as general quasismooth wellformed hypersurfaces of trivial canonical class in a 4-dimensional weighted projective space.}

\bibliographystyle{alpha}
\bibliography{biblio}

\end{document}